\documentclass[11pt]{book}
\usepackage[utf8]{inputenc}
\usepackage{amsfonts, amsmath, amssymb, amsthm}  
\usepackage{times}
\usepackage[english]{babel}
\usepackage{graphicx}
\usepackage{biblatex}
\usepackage{comment}
\usepackage{csquotes}
\DeclareGraphicsExtensions{.pdf,.png,.jpg}
\usepackage{mathtools}
\usepackage{mathrsfs}
\newtheorem{theorem}{Theorem}[section]

\newtheorem{definition}{Definition}[section]

\begin{document}
\pagestyle{empty}

\begin{center}

\vspace{1cm}

{\Huge         Projective Geometry and PDE Prolongation}

\vspace{35mm} 

\includegraphics[width=2cm]{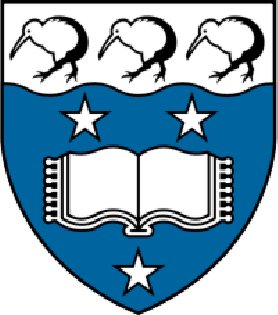}

 \vspace{45mm}

{\Large       Jake McNaughton}

	\vspace{1ex}

Department of Mathematics

The University of Auckland

	\vspace{5ex}

Supervisor:            Professor Rod Gover

	\vspace{30mm}

A dissertation  submitted in partial fulfillment of the requirements for the degree of BSc(Hons)  in Mathematics, The University of Auckland, 2021.

\end{center}

  \newpage

\chapter*{Abstract}    
In this dissertation we study basic local differential geometry, projective differential geometry, and prolongations of overdetermined geometric partial differential equations. It is simple to prolong an $n$-th order linear ordinary differential equation into $n$ first order equations. For partial differential equations there is a related process but it is far more subtle and complex. Considerable work has been done in this area but much of this is very abstract and there are many open problems even at a relatively elementary level.  We introduce the reader to differential geometry and tractor calculus before recovering the projective tractor and cotractor connections via the prolongation of appropriate partial differential equations. Following this, we study prolongation of other projectively invariant equations, a particular focus is an equation known as the projective metrisability equation.
\setcounter{page}{1}
\pagestyle{headings}
\pagenumbering{roman}

\addcontentsline{toc}{chapter}{Abstract}
\setcounter{secnumdepth}{3} 
\setcounter{tocdepth}{3}
\tableofcontents            


\chapter*{Introduction}
\addcontentsline{toc}{chapter}{Introduction}
\markboth{INTRODUCTION}{INTRODUCTION}


%

Between 1925 and 1926, Tracey Thomas developed a calculus for projective and conformal manifolds \cite{Thomas588}. Thomas's work on the subject was subsequently forgotten until 1994 when Bailey, Eastwood, and Gover rediscovered and expanded it in their paper titled `Fall 1994
Thomas's Structure Bundle for Conformal, Projective and Related Structures' \cite{BEGogThomassBundle}. After its rediscovery, Thomas's methods of calculus were named ``tractor'' calculus, a portmanteau of Tracey Thomas and vector, to acknowledge Thomas \cite{BEGogThomassBundle}. Projective differential geometry started being developed in the 1920s \cite{EastwoodNotesonPDG} with large contributions made by Tracey Thomas and \'Elie Cartan. It has since fallen into obscurity. Eastwood's notes \cite{EastwoodNotesonPDG} give a great overview of the topic. In undergraduate courses on differential equations it is common to reduce an $n$-th order ordinary differential equation to a system of $n$ linear first order equations, this process is called prolongation. Prolongation of PDEs is much more subtle and complicated \cite{BEAGprolongation}.

\vspace{0.2cm}
\noindent
Projective tractor calculus can be linked to the prolongation of a particular overdetermined differential equation. The prolongation yields a distinguished connection on a distinguished bundle called the tractor bundle. From this, we can extract information about the solution space of the equation that has been prolonged. Prolongation can be applied to a variety of equations - specifically overdetermined geometric partial differential equations - and provides information about the solution spaces. 

\vspace{0.2cm}
\noindent
In this dissertation, we will study prolongation. This will involve describing the prolongation process with many examples included to familiarise the reader with it. Using tractor calculus we will see how distinguished connections on higher order bundles can be found without performing any further prolongations, and how dual connections can be calculated using the Leibniz rule. This dissertation bridges gaps  between aspects of introductory differential geometry texts and literature introducing projective geometry and tractor calculus, which tends to be aimed at experts in the field.

\vspace{0.2cm}
\noindent
We begin in Chapter 1 by introducing core concepts in differential geometry, focusing on working towards the definitions of Riemannian geometry and projective differential geometry. We will provide all necessary knowledge of differential geometry, beginning with the definition of a manifold. 

\pagebreak

\noindent
In Chapter 2, we introduce tractor calculus. The key object here is the tractor connection and we will be working in the projective case with the projective tractor connection. After the tractor connection and its dual (the cotractor connection) are defined, we give the prolongations leading to them.

\vspace{0.25cm}
\noindent
In Chapter 3, we study prolongation of projectively invariant equations. We prolong the projective metrisability equation

\begin{equation}
    \text{trace-free part of }({\nabla_{a}t^{bc}})=0.
\end{equation}

\noindent Non-degenerate solutions to this equation are equivalent to metrics whose Levi-Civita connection lies in the projective class. The connection obtained by prolonging the metrisability equation is the same as the tractor connection acting on that bundle if and only if the metric is Einstein. This is an important result as Einstein metrics have a special role in projective differential geometry \cite{GoverEinsteinGeodesics}.

\vspace{0.3cm}
\noindent
This dissertation has been designed to be accessible to those not well-versed in differential geometry. We provide a self-contained work that brings together existing knowledge on the subject.

\chapter{Preliminaries}
\pagenumbering{arabic}
In this chapter we give an overview of core concepts in differential geometry, building up to a definition of Riemannian geometry. We then proceed to a brief introduction to projective differential geometry, mostly based on `Notes on Projective Differential Geometry' by Michael Eastwood \cite{EastwoodNotesonPDG}, which the reader is pointed to for further details.
\section{Notation}
We will predominantly use Roger Penrose's abstract index notation interspersed with index-free notation when it is more convenient. This allows us to represent tensors in a concise and clear, way with an upper index representing a contravariant component and a lower index representing a covariant component. An $\begin{psmallmatrix}
n\\m
\end{psmallmatrix}$ tensor over a vector space $V^{a}$, which will be defined later in this chapter, can be viewed as a multilinear map from $m$ direct products of the vector space, $V^{a}$, and $n$ direct products of the dual vector space, $V_{a}$ or $V^{\ast}$ into $\mathbb{R}$, and hence is represented with $n$ upper indices and $m$ lower indices. For example, $T_{ab}{}^{c}{}_{d}{}^{ef}:V^{\ast}\times V^{\ast}\times V\times V^{\ast}\times V\times V\to \mathbb{R}$. An $\begin{psmallmatrix}
n\\m
\end{psmallmatrix}$
 tensor field maps into $\mathbb{R}$ at each point on the manifold and therefore maps from the tensor bundle into real-valued functions. When an index appears as an upper index and lower index it represent a contraction, as in the dual pairing of the vector and covector components in question.

We use $TM$, but also $\mathcal{E}^{a}$, to denote the tangent bundle of a manifold $M$, and similarly $T^{\ast}M$ or $\mathcal{E}_{a}$ to denote the cotangent bundle of a manifold $M$. In section 2.1 we will introduce density bundles, weighted sections of these bundles with weight $k$ are denoted $\mathcal{E}^{a}(k)$. For indices of the tangent bundle and its element we will use lower case latin letters, and for the indices of the tractor bundle and its elements we will use upper case latin letters. We denote the tractor bundle $\mathcal{E}^{A}$ and the tractor connection $\nabla^{\mathcal{T}}_{a}$. $X^{A}, Y_{A}, Z_{A}^{a}, W^{A}_{a}$ denote splitting maps for the standard projective tractor connection when appropriate, however we will also use $X, Y, Z$ and $W$ as arbitrary variables when necessary.

We use $[A,B]=AB-BA$ to denote the Lie bracket of vector fields $A$ and $B$. We denote sections of a vector bundle $V$ as $\Gamma(V)$, for example, sections of the tangent bundle are denoted $\Gamma(TM)$ or $\Gamma(\mathcal{E}^{a})$.

\section{Manifolds, the Tangent Bundle, and Connections}
This section provides a concise explanation of the most basic object in differential geometry - the manifold. We begin with the definition of a manifold and then build up some structure in the form of vectors on our manifold.
\begin{definition}[Manifold]
A $C^{\infty}$ $n$-dimensional manifold $M$ is a set $M$ together with a $C^{\infty}$ atlas $\{U_{\alpha},\phi_{\alpha}\}$, that is to say a collection of charts $(U_{\alpha},\phi_{\alpha})$ where the $U_{\alpha}$ are subsets of $M$ and the $\phi_{\alpha}$ are one-to-one maps from the corresponding $U_{\alpha}$ to open sets of $\mathbb{R}^{n}$ such that

\begin{itemize}
\item[(i)] the $U_{\alpha}$ cover $M$, i.e. $M=\bigcup_{\alpha}U_{\alpha}$,
\item[(ii)] if $U_{\alpha}\cap U_{\beta}$ is non-empty, then the map $$\phi_{\alpha}\circ \phi_{\beta}^{-1}:\phi_{\beta}(U_{\alpha}\cap U_{\beta})\to \phi_{\alpha}(U_{\alpha}\cap U_{\beta})$$ is a $C^{\infty}$ map of an open subset of $\mathbb{R}^{n}$ to an open subset of $\mathbb{R}^{n}$.
\end{itemize}
\end{definition}
\noindent Equipped with this basic definition \cite{hawking_ellis_1973} we will proceed  by defining a vector bundle \cite{O’Neill}. We will be very interested in a specific vector bundle - the tangent bundle.

\begin{definition}[Vector Bundle]
A real rank $k$ vector bundle $(E,\pi)$ over a manifold $M$ - also called the base space - consists
of a manifold $E$ and a smooth map $\pi : E \to M$ such that 

(1) for each $\pi^{-1}(p),p\in M$, is a $k$-dimensional vector space;

(2) for each $p \in M$ there is a neighbourhood
$\mathscr{U}$ of $p$ in $M$ and a diffeomorphism
$$\phi:\mathscr{U}\times \mathbb{R}^{k}\to\pi^{-1}(\mathscr{U})\subset E$$
such that for each $q\in\mathscr{U}$, the map $v\to\phi(q,v)$ is a linear isomorphism from $\mathbb{R}^{k}$ onto $\pi^{-1}(q)$.
\end{definition}

\begin{definition}[Section of a Vector Bundle]
A section of a vector bundle $V$ is a smooth, injective map $X : B\to E$ such that $\pi\circ X=Id_{M}$
\end{definition}

\noindent Now we are going to look at the tangent bundle. At each point $p$ on a manifold, $M$, we can define its tangent space which gathers all tangent vectors at that point. Hence we will proceed by first defining a tangent vector \cite{hawking_ellis_1973}.

\begin{definition}[Tangent Vector]
Let $k\geq 1$. A $C^{k}$ curve, $\gamma(t)$, in $M$ is a $C^{k}$ map from $\mathbb{R}$ into $M$. The vector, $(\tfrac{\partial }{\partial t})_{\gamma}\vert_{t_{0}}$, tangent to the $C^{1}$ curve $\gamma(t)$ at the point $\gamma(t_{0})$ is the operator which maps each $C^{1}$ function $f$ into the number $(\tfrac{\partial f}{\partial t})_{\gamma}\vert_{t_{0}}$. All vectors that come about in this way are the tangent vectors to $M$ at $t_{0}$.
 
\end{definition}

\noindent This coincides with the usual understanding of a tangent vector to a curve or surface. Now with this definition we can define tangent spaces \cite{Gallot} on a manifold which consist of vector spaces containing all of the tangent vectors at each point.

\begin{definition}[Tangent Space]
The tangent space to $M$ at $p$, denoted $T_{p}M$, is the set of all tangent vectors to $M$ at $p$. The tangent space at each point on an $n$-dimensional manifold naturally gets the structure of an $n$-dimensional manifold.
\end{definition}

\noindent The tangent bundle is a global structure on the manifold comprised of all the tangent spaces on $M$. 
\begin{definition}[Tangent Bundle]
If $M$ is a differentiable manifold, the tangent bundle $TM$ is a vector bundle such that for all $p\in M$, $\pi^{-1}(p)=T_{p}M$, is the tangent space at $p$, where $\pi$ is the bundle map. This is equipped with a smooth structure characterised by the property that coordinate vector fields are locally smooth sections.

\end{definition}

\noindent The dual vector space to $T_{p}M$ is called the cotangent space, denoted $T_{p}^{\ast}M$. Elements of $T_{p}^{\ast}M$ are called cotangent vectors. The disjoint union of all cotangent spaces over $M$ equipped with a smooth structure - similar to what we have for the tangent bundle - is called the cotangent bundle and denoted $T^{\ast}M$.

\subsection{Sections of the Tangent Bundle}
A section of the tangent bundle is a smooth, injective map $X:M\to TM$ such that $\pi\circ X=Id_{M}$ . Therefore, a section smoothly identifies one vector in each tangent space - that is $X(p)\in T_{p}M$. We denote the set of sections on the tangent bundle of $M$ by $\Gamma(TM)$, and this is set of all vector fields on $M$.  

We continue by introducing tensors which generalise vectors to higher dimensions.

\pagebreak

\section{Tensors}
\subsection{Tensors as Multilinear Maps}
\begin{definition}[Tensor]
Let $p$ and $q$ be integers, not both zero, and $V$ be a vector space. A multilinear function:
$$T:\underbrace{V^{\ast}\times\dots \times V^{\ast}}_{p\text{ copies}}\times \underbrace{V\times \dots\times V}_{q\text{ copies}}\to \mathbb{R}$$
is called a tensor of type $\left(\begin{smallmatrix}p\\q\\\end{smallmatrix}\right)$ over $V$. A tensor of type $\left(\begin{smallmatrix}0\\0\\\end{smallmatrix}\right)$ is defined to be an element of $\mathbb{R}$. \cite{shlomo}

\end{definition}

\subsection{Tensor Fields}

A tensor field assigns a tensor to each point in a space. As we are building up an understanding of manifolds, naturally this is the setting in which we will be considering tensor fields. A tensor field, $T$, on a manifold $M$, assigns a tensor to each point of $M$. We should note that we are primarily concerned with tensor fields over the tangent bundle, so the tensor assigned to $p\in M$ will be a tensor over the vector space $T_{p}M$.

\subsubsection{An Example: The Metric Tensor}
An example of a tensor field is the metric tensor, which is a type $\left(\begin{smallmatrix}0\\2\\\end{smallmatrix}\right)$ tensor field over $TM$, and gives a notion of distance on $M$. Despite commonly just being referred to as a metric, a metric tensor is a distinct but related concept to the distance function of a metric space.

A metric tensor that is positive-definite is an inner product and generalises the dot product to a manifolds, allowing the length of curves on the manifold to be calculated. Thus, it defines a metric on the manifold where the distance between two points on the manifold is the infimum of the lengths of all the curves which join the points.

\subsection{Another Perspective on Tensors}
An alternative view of tensors, which is common in physics, is as multidimensional arrays with certain properties. Tensors are a generalisation of vectors and we can write an $n-$dimensional vector as an $n\times 1$ matrix in a certain basis, we can do something similar with a rank $q$ tensor. In a certain basis, a rank $q$ tensor can be thought as a $q$-dimensional hypercube which follows particular transformation rules under a change of basis.

\subsection{Differentiating Tensors}
To develop calculus methods that work for the concepts we have been introducing we need a way of comparing vectors in different tangent spaces. This would give a way to differentiate vectors over a manifold and also a way of doing this for higher rank tensors. The device we use for this is called a connection. It will allow us to compare nearby tangent spaces. A connection at a point $p\in M$ is a rule which assigns to each vector field $X$ at $p$ a differential operator $\nabla_{X}$ which maps a vector field $Y$ into a vector field $\nabla_{X}Y$.

It is common to use a connection defined on the tangent bundle of a manifold, however later we will be working with connections on a vector bundle called the tractor bundle (of a manifold, $M$) which is of a different rank to the tangent bundle. Hence, we will begin with the definition of a connection on an arbitrary vector bundle and then give the definition of a connection on the tangent bundle \cite{Gallot}.

\begin{definition}[Connection on a Vector Bundle]
	Let $M$ be a smooth manifold and $E$ be a vector bundle over $M$. A connection on $E$ is a bilinear map from $\Gamma(TM)\times \Gamma(E)$ to $\Gamma(E)$ such that for $X\in \Gamma(TM), Y\in\Gamma(E)$ and $f\in C^{\infty}(M)$, we have
    $$\nabla_{fX}Y=f\nabla_{X}Y\quad \text{and} \quad \nabla_{X}fY=(Xf)Y+f\nabla_{X}Y$$
\end{definition}

\noindent Every connection on a vector bundle is an affine connection. We can think of $\nabla_{X}Y$ as the covariant derivative of $Y$ in the direction of $X$ \cite{moore2006lectures}. A connection on the tangent is simply a connection on the vector bundle $E=TM$. For completeness we provide a definition.

\begin{definition}[Connection on the Tangent Bundle] A linear connection on the tangent bundle $TM$ is a map $$\nabla:\Gamma(TM)\times \Gamma(TM)\to\Gamma(TM)$$
Such that for $X,Y,Z,W\in\Gamma(TM)$, $\alpha,\beta\in\mathbb{R}$ and $f,g\in C^{\infty}(M)$, $\nabla$ satisfies

\begin{itemize}
    \item[(i)] $\nabla_{fX+gW}Y=f\nabla_{X} Y+g\nabla_{W} Y$
    \item[(ii)] $\nabla_{X}(\alpha Y+\beta Z)=\alpha\nabla_{X}Y+\beta\nabla_{X}Z$
    \item[(iii)] $\nabla_{X}(fY)=X(f)Y+f\nabla_{X} Y$
\end{itemize}

\end{definition}
\noindent Property $(iii)$ is called the Leibniz rule. 

\subsubsection{Duality and Connections on Tensor Products}

The Leibniz rule holds for arbitrary tensor fields $S$ and $T$ in the following way, $\nabla_{a}(S\otimes T)=\nabla_{a}S\otimes T+S\otimes \nabla_{a}T$. This allows us to find how the connection acts on elements of the dual bundle by applying the Leibniz rule to a covariant vector contracted with a contravariant vector as such, $\nabla_{a}(\mu_{b}\nu^{b})=\mu_{b}\nabla_{a}\nu^{b}+\nu^{b}\nabla^{\ast}_{a}\mu_{b}$. Also, $\nabla f=df$ for any function $f$, so we get $d_{a}(\mu_{b}\nu^{b})=\mu_{b}\nabla_{a}\nu^{b}+\nu^{b}\nabla^{\ast}_{a}\mu_{b}$ which defines the connection $\nabla^{\ast}$ on $T^{\ast}M$. We will now denote both these connections as $\nabla$.

Similarly, we can use the Leibniz rule to see how the connection acts on higher order tensors. In this way, we can find how the connection acts on all tensors over the tangent and cotangent bundle of the manifold.

\section{Properties of a Connection on the Tangent Bundle}

There are two important tensors that an affine connection $\nabla$ defines, namely \textbf{torsion} and \textbf{Riemannian curvature} which shall be defined here.
\begin{definition}[Torsion] The torsion of an affine connection, $\nabla$, on a manifold $M$ is the type $\left(\begin{smallmatrix}1\\2\\\end{smallmatrix}\right)$ tensor which is given by $T^{\nabla}(X,Y)=\nabla_{X}Y-\nabla_{Y}X-[X,Y]$ where $X$ and $Y$ are sections of the tangent bundle of $M$.

\end{definition}
\noindent A \textbf{torsion-free} connection is one whose torsion vanishes for each $X$ and $Y$ in $\Gamma(TM)$, we will be working extensively with torsion-free connections.

\begin{definition}[Riemannian Curvature]The Riemannian curvature of a connection, $\nabla$ acting on a vector bundle $V$is the type $\left(\begin{smallmatrix}1\\3\\\end{smallmatrix}\right)$ tensor which is given by $R^{\nabla}(U,V)=[\nabla_{U},\nabla_{V}]-\nabla_{[U,V]}$ where $U$ and $V$ are vector fields on $M$.
\end{definition}
\noindent A \textbf{flat} connection is one whose Riemannian curvature vanishes everywhere on the manifold. We can form another curvature from the Riemannian curvature, namely the Ricci curvature.
\begin{definition}[Ricci Curvature] 
The Ricci curvature is a type $\left(\begin{smallmatrix}0\\2\\\end{smallmatrix}\right)$ tensor defined by contracting over the first and third indices of the Riemannian curvature and is denoted $R_{ab}=R_{ca}{}^{c}{}_{b}$.
\end{definition}

\subsection{Symmetries of the Riemannian Curvature}
The Riemannian curvature tensor has many symmetries. Firstly, it is skew-symmetric in its first two indices:
$$R_{ab}{}^{c}{}_{d}=-R_{ba}{}^{c}{}_{d}.$$
Additionally, it has the symmetry

$$R_{ab}{}^{c}{}_{d}+R_{da}{}^{c}{}_{b}+R_{bd}{}^{c}{}_{a}=0$$
which is known as the first Bianchi identity. Finally, it also obeys the second Bianchi identity,

$$\nabla_{e} R_{ab}{}^{c}{}_{d}+\nabla_{b} R_{ea}{}^{c}{}_{d}+\nabla_{a} R_{be}{}^{c}{}_{d}=0.$$
\textbf{Remark:} The second Bianchi identity holds for a connection over any vector bundle as it follows from the Leibniz rule and the Jacobi Identity, however for vector bundles other than the tangent bundle it will be in a different form to the above.


\section{Riemannian Geometry}
Riemannian geometry is the study of Riemannian manifolds, smooth manifolds with a Riemannian metric.
\begin{definition}[Riemannian Manifold]
A Riemannian manifold $(M, g)$ is a manifold, $M$, endowed with a positive-definite metric $g$.
\end{definition}

\noindent We discussed metrics earlier but will give a precise definition here.
\begin{definition}[Metric Tensor]
A metric tensor, $g$, is a type $\begin{pmatrix}
0\\2
\end{pmatrix}$ tensor field on $M$ which at every point is symmetric and positive definite.
\end{definition}

\noindent \textbf{Remark:} The metric $g$ at $p\in M$ assigns a magnitude $|g(X,X)|^{\frac{1}{2}}$ to each $X\in T_{p}M$ and defines the angle $$\frac{g(X,Y)}{|g(X,X)\cdot g(Y,Y)|^{\frac{1}{2}}}$$
between $X,Y\in T_{p}M.$ In this way, the metric is a generalisation of the dot product in Euclidean space.

Riemannian Geometry is the study of Riemannian manifolds. If we relax the condition of the metric being positive definite then we have a Psuedo-Riemannian manifold. On every Riemannian or Psuedo-Riemannian manifold there exists a distinguished connection called the \textbf{Levi-Civita Connection} which is distinguished by the metric. We will give the definition of this connection.

\begin{definition}[Levi-Civita Connection] The Levi-Civita Connection is the unique affine connection that exists on any Riemannian manifold $(M,g)$ that is torsion-free and is compatible with the metric (This means that $Xg(Y,Z)=g(\nabla_{X}Y,Z)+g(Y,\nabla_{X}Z)$ for any $X,Y,Z\in\Gamma(TM)$. In abstract index notation this is denoted $\nabla_{k} g_{ij}=0$).
\end{definition}

\noindent The existence and uniqueness of such a connection is proved in \ref{appendixB}.

\subsection{Raising and Lowering Indices with the Metric Tensor}
The metric and inverse metric can be contracted with tensors to lower and raise indices respectively. The metric lowers indices in the following way, $g_{ab}X^{a}=X_{b}$, and the inverse metric raises indices by $g^{ab}X_{a}=X^{b}$. This allows us to define another type of curvature.

\begin{definition}
[Scalar Curvature]
The Scalar curvature is a type $\left(\begin{smallmatrix}0\\0\\\end{smallmatrix}\right)$ tensor or scalar defined by using the inverse metric tensor to raise the index in the following way $g^{ab}R_{ac}=R^{b}{}_{c}$ and then contracting over the two indices by taking the trace $\delta_{b}^{c}R^{b}{}_{c}=R^{b}{}_{b}$ It is denoted $R=Ric^{b}{}_{b}$.
\end{definition}

\section{Projective Differential Geometry} 
Projective differential geometry provides a simple setting in which to study overdetermined systems of partial differential equations which we will do in subsequent chapters. It is the study of manifolds with an additional projective structure. 
\subsection{Geodesics}
If we equip a manifold with a connection, then geodesics are the curves on the manifold which represent the shortest paths between points. In flat space these would be straight lines however, for example on a sphere they are the great circles. In projective geometry we are interested in connections which all have the same unparameterised geodesic curves.

\begin{theorem}
Two torsion-free connections $\nabla$ and $\overline{\nabla}$ have the same geodesics as unparameterised curves if and only if for any section of one-forms, $\omega_{a}$,
$$\overline{\nabla}_{a} \omega_{b}=\nabla_{a}\omega_{b}-\Upsilon_{a}\omega_{b}-\Upsilon_{b}\omega_{a}$$ for some one-form $\Upsilon_{a}$.
\end{theorem}
\noindent We call two torsion-free connections \textbf{projectively equivalent} if they have the same geodesics when expressed as unparameterised curves. 

\begin{definition}[Projective Structure] A projective structure on $M$ is an equivalence class of torsion-free projectively equivalent connections.

\end{definition}
\noindent In projective geometry we study manifolds equipped with a projective structure.

\nocite{*}

\subsection{Projective Invariance}
Something which is projectively invariant is well-defined on a manifold with just a projective structure on it. The subtlety here is that often projectively invariant objects are expressed in terms of a connection in the projective class. Such a projectively invariant object will remain the same when the connection is switched to any other connection within the projective class.

We will be working exclusively with special affine connections. For special affine connections, the projective Schouten tensor is symmetric. The Schouten tensor is a projectively invariant tensor given by $P_{ab}=\frac{1}{n-1}R_{ab}$. We will also be using the Weyl tensor which is the trace free component of the Riemannian curvature $W_{ab}{}^{c}{}_{d}=R_{ab}{}^{c}{}_{d}-\delta_{a}^{c}P_{bd}+\delta_{b}^{c}P_{ad}$, this is also projectively invariant. Finally, the Cotton tensor $C_{abc}=\nabla_{a}P_{bc}-\nabla_{b}P_{ac}$ is also projectively invariant. It is related to the Weyl tensor by $(n-2)C_{dab}=\nabla_{c}W_{ab}{}^{c}{}_{d}$.

\chapter{Tractor Calculus}
A projective structure does not have a distinguished connection on the tangent bundle such as in the case of Riemannian geometry with the Levi-Civita connection. However, it turns out that there does exist a distinguished connection (the tractor connection) on a higher rank bundle (the tractor bundle) which has rank $n+1$.

We begin by defining the tractor bundle. Then we see how the tractor connection acts on elements of this bundle. With this, the Leibniz rule and properties of tensor powers define a dual connection acting on the so called `cotractor bundle' and also define the way which these connections act on tensors of higher rank over these two bundles.

\section{Weighted Bundles}
We have already seen the prototypical example of a vector bundle: the tangent bundle. In Riemannian geometry this is the bundle over which tensor fields are defined and what the connections act on. Smooth sections of this bundle are vector fields. Now we will be looking at higher rank bundles.

\subsection{Density Bundles}
As introduced in the previous chapter, $TM$ is the tangent bundle of the manifold $M$. The notation $\Lambda^{k}TM$ represents the $k$-th exterior power of the tangent bundle, when $k=n$ we have $\Lambda^{n}TM$ which is called the \textbf{top exterior power}. Then $(\Lambda^{k}TM)^{2}$ is an oriented line bundle denoted by $\mathcal{K}$ and called the density bundle of weight $2n+2$.

Since $\mathcal{K}$ is oriented there is a canonical way to take oriented roots, so we define the density bundle of weight $w$ to be $\mathcal{E}(w)=\mathcal{K}^{\frac{w}{2n+2}}$. For example, the density bundle of weight $1$, $\mathcal{E}(1)$, is $\mathcal{K}^{\frac{1}{2n+2}}$. Now, if we take the tensor product of a vector bundle, $\mathcal{B}$, with the density bundle of weight $w$ it is common to use the shortened notation $\mathcal{B}(w)\coloneqq \mathcal{B}\otimes \mathcal{E}(w)$. 

A line bundle is a vector bundle of dimension one. Any density bundle is naturally a line bundle, and sections of the density bundle of weight $w$ are called densities of weight $w$.
\section{The Projective Tractor Bundle}

We consider a specific vector bundle of rank $n+1$ as this is the one where we can find a distinguished connection on for a projective manifold $(M,p)$.

\begin{definition}[Cotractor Bundle]
For each choice of connection in the projective class, $\mathcal{E}_{A}$ is identified as the direct sum $\mathcal{E}(1) \oplus \mathcal{E}_{a}(1)$. When the connection is changed to another in the projective class, this splitting changes according to $$\overline{\begin{pmatrix} \sigma\\ \mu_{a}\end{pmatrix}}=\begin{pmatrix} \sigma \\ \mu_{a}+\Upsilon_{a}\sigma\end{pmatrix}.$$ Where $\mathcal{E}(1)$ denotes a root of the top exterior power and gives the functions of weight 1. $\mathcal{E}_{a}(1)$ is the bundle of 1-forms of weight 1 obtained by taking the tensor product of the 1-forms with $\mathcal{E}(1)$. $\Upsilon_{a}$ is the one form from Theorem 1.6.1. \cite{EastwoodNotesonPDG}

\end{definition}

\noindent The dual bundle to this is the tractor connection and is identified as the direct sum $\mathcal{E}^{A}=\mathcal{E}(-1)\oplus\mathcal{E}^{a}(-1)$. Now we have this bundle which we will show has a projectively invariant distinguished connection.

We will make choices of connection in the projective class, thus obtaining a splitting, without further comment. Given a choice of connection we write $Z_{A}^{a}:\mathcal{E}_{a}(1)\to\mathcal{E}_{A}$, $X^{A}:\mathcal{E}_{A}\to\mathcal{E}(1)$, $W_{a}^{A}:\mathcal{E}_{A}\to\mathcal{E}_{a}(1)$ and $Y_{A}:\mathcal{E}(1)\to\mathcal{E}_{A}$ for the corresponding bundle maps which give the splitting \cite{NeusserSasaki}. Then elements of $\mathcal{E}_{A}$ are of the form $\begin{psmallmatrix}
\sigma\\\mu_{b}
\end{psmallmatrix}=Z_{A}^{a}\mu_{a}+Y_{A}\sigma$, and elements of $\mathcal{E}^{A}$ are of the form $\begin{psmallmatrix}
\nu^{b}\\t
\end{psmallmatrix}=W_{a}^{A}\nu^{a}+X^{A}t$.
\section{The Projective Tractor Connection}
The projective tractor connection is the distinguished connection which acts on the bundle given above. In this section we give how the projective tractor connection acts on elements of the tractor bundle and dual tractor bundle. The projective tractor connection will be denoted as $\nabla_{a}^{\mathcal{T}}$, it acts on an element of the tractor bundle in the following way,

\begin{equation}\label{dualprojtrac}
    \nabla_{a}^{\mathcal{T}}\begin{pmatrix}
    \mu^{b}\\t
    \end{pmatrix}=\begin{pmatrix}
    \nabla_{a}\mu^{b}+t\delta_{a}^{b}\\
    \nabla_{a}t-P_{ab}\mu^{b}
    \end{pmatrix}.
\end{equation}

\noindent Similarly, we denote the dual tractor connection by $\nabla_{a}^{\mathcal{T}\ast}$ which acts on the cotractor bundle in the following way,
\begin{equation}\label{projtrac}
    \nabla_{a}^{\mathcal{T}\ast}\begin{pmatrix}
    \sigma\\\nu_{b}
    \end{pmatrix}=\begin{pmatrix}
    \nabla_{a}\sigma-\nu_{a}\\
    \nabla_{a}\nu_{b}+P_{ab}\sigma
    \end{pmatrix}.
\end{equation}

\noindent These are both connections on their respective bundles. 

\subsection{Curvature of the Tractor Connection}
The curvature of the dual tractor connection is given by

$$\nabla_{a}\nabla_{b}-\nabla_{b}\nabla_{a}\begin{pmatrix}\sigma\\\mu_{c}\end{pmatrix}=\begin{pmatrix}0\\-W_{ab}{}^{d}{}_{c}\mu_{d}+C_{abc}\sigma\end{pmatrix},$$

\noindent\cite{EastwoodNotesonPDG}. On the other hand,

$$\nabla_{a}\nabla_{b}-\nabla_{b}\nabla_{a}\begin{pmatrix}\nu^{c}\\\rho\end{pmatrix}=\begin{pmatrix}W_{ab}{}^{c}{}_{d}\nu^{d}\\-C_{ab}d\nu^{d}\end{pmatrix}$$

\noindent gives the curvature of the tractor connection \cite{BEGogThomassBundle}.
\noindent is the curvature of the tractor connection

\section{Prolongation and the Tractor Connection}
The projective tractor bundle and connection arise out of a process called prolongation. The prolongation of a specific projectively invariant equation yields the projective tractor connection.  Connections arising from other equations through prolongation are also called tractor connections.

\subsection{Prolongation Leading to the Tractor Connection}
We wish to prolong $\text{tf}(\nabla_{a}\nu^{c})=0$, as this is the equation that yields the projective tractor connection. The equation states that the trace-free part of $\nabla_{a}\nu^{c}$ vanishes, so we begin by finding an expression for the trace-free part. Therefore, we let $\delta_{a}^{c}\mu$ be the trace part of $\nabla_{a}\nu^{c}$. Hence,
\begin{equation}\label{projtractoreq1}
    \nabla_{a}\nu^{c}-\delta_{a}^{c}\mu=0,
\end{equation}
since $\text{tf}(\nabla_{a}\nu^{c})=\nabla_{a}\nu^{c}-\delta_{a}^{c}\mu$. To find out what $\mu$ is we can take the trace over $a$ and $c$, $\delta_{c}^{a}\nabla_{a}\nu^{c}-\delta_{c}^{a}\delta_{a}^{c}\mu=0$. Since $a$ and $c$ are contracted and run from $1$ through to $n$, contributing $1$ to the value each time, hence $\delta_{c}^{a}\delta_{a}^{c}=\sum^{n}_{1} 1=n$. Thus, $\nabla_{c}\nu^{c}-n\mu=0$ which means $\mu=\frac{1}{n}\nabla_{c}\nu^{c}$. Now we proceed with the prolongation by taking the derivative using the connection,

$$\nabla_{b}\nabla_{a}\nu^{c}-\nabla_{b}(\delta_{a}^{c}\mu)=0.$$

\noindent Our aim is to find an equation with $\nabla_{a}\mu$ in terms of known terms. To do this, we wish to swap the order of the covariant derivatives $\nabla_{a}$ and $\nabla_{b}$. 
To swap this we must introduce a curvature term , since $\nabla_{b}\nabla_{a}\nu^{c}=\nabla_{a}\nabla_{b}\nu^{c}-R_{ab}{}^{c}{}_{d}\nu^{d}$, so as we have it,
$$\nabla_{a}\nabla_{b}\nu^{c}-R_{ab}{}^{c}{}_{d}\nu^{d}-\nabla_{b}(\delta_{a}^{c}\mu)=0.$$
Now we want to trace over $b$ and $c$ so that the first term becomes a multiple of $\nabla_{a}\mu$, this gives $\delta_{c}^{b}(\nabla_{a}\nabla_{b}\nu^{c}-R_{ab}{}^{c}{}_{d}\nu^{d}-\nabla_{b}(\delta_{a}^{c}\mu))=0$, and now we rearrange.

$$\nabla_{a}\nabla_{c}\nu^{c}-R_{ac}{}^{c}{}_{d}\nu^{d}-\delta_{a}^{c}\nabla_{c}\mu=0$$

\noindent Using the fact that $\nabla_{c}\nu^{c}=n\mu$ and symmetries of the Riemannian curvature yields,

$$n\nabla_{a}\mu+R_{ca}{}^{c}{}_{d}\nu^{d}-\nabla_{a}\mu=0$$

\noindent Next we rearrange and replace the contracted Riemannian curvature with the Ricci curvature to reach the final equation.

$$(n-1)\nabla_{a}\mu-R_{ad}\nu^{d}=0$$

\begin{equation}\label{projtractoreq2}
    \implies \nabla_{a}\mu+\frac{R_{ad}\nu^{d}}{n-1}=0
\end{equation}

\noindent Equations \ref{projtractoreq1} and \ref{projtractoreq2}  now form a closed system so can be brought together to give the projective tractor connection,
\begin{equation}\label{projtractorconnection}
    D_{a}\begin{pmatrix}
    \nu^{c}\\
    \mu
    \end{pmatrix}=\begin{pmatrix}
    \nabla_{a}\nu^{c}-\delta_{a}^{c}\mu\\
    \nabla_{a}\mu+ \frac{R_{ad}\nu^{d}}{n-1}
    \end{pmatrix}.
   \end{equation}
   

\subsection{Prolongation Leading to the Dual Tractor Connection}
To familiarise the reader with the prolongation before studying it more extensively, we will prolong the equation that leads to the dual projective tractor connection. In the flat case, this equation is $\nabla_{a}\nabla_{b}\sigma=0$, however in the curved case there is another term. In the curved case, the equation needing to be prolonged is 
\begin{equation}\label{dualeq0}
    \nabla_{a}\nabla_{b}\sigma+P_{ab}\sigma=0.
\end{equation}
where $P_{ab}$ denotes the projective Schouten tensor.
In the previous prolongation there was an obvious place to start, namely, finding the trace-free part. In contrast, the starting place here is less obvious. Prolongation requires defining new variables. Since we want our equations to be in terms of the variables and their first derivatives we need to deal with the $\nabla_{a}\nabla_{b}\sigma$ term so that the second derivative does not appear. To do this we define $\nu_{b}\coloneqq\nabla_{b}\sigma$. This gives us our first equation, namely,
\begin{equation}\label{dualeq1}
    \nabla_{b}\sigma-\nu_{b}=0.
\end{equation}

\noindent Now we have an equation involving $\nabla_{b}\sigma$, we need one involving $\nabla_{b}\nu_{c}$. We get this immediately by substituting $\nu_{b}=\nabla_{b}\sigma$ to yield

\begin{equation}\label{dualeq2}
    \nabla_{a}\nu_{b}+P_{ab}\sigma=0.
\end{equation}

\noindent This completes the prolongation, and the equations \ref{dualeq1} and \ref{dualeq2} form the prolonged system, giving the standard tractor connection as

\begin{equation*}
   \nabla_{a}\begin{pmatrix}
    \sigma\\\nu_{b}
    \end{pmatrix}=\begin{pmatrix}
    \nabla_{a}\sigma-\nu_{a}\\
    \nabla_{a}\nu_{b}+P_{ab}\sigma
    \end{pmatrix}.
\end{equation*}

\subsection{Prolongation in its Own Right}
The rest of this dissertation will be devoted to studying the method of prolongation. We will see how this can give information about the equations which are prolonged and how it can indicate links between equations. We will give a lot of attention to duality and how that relates to prolonging equations and the resultant connections. How the prolonged system gives information about solutions to the original equation will also be looked at. Finally, we will look at how the tractor connection can be used to assist in this process. 
\newpage 
 \chapter{Prolongation}
In Chapter 2, the prolongations of two equations were done to obtain the projective tractor connection and its dual. This indicates that there is some sort of underlying link between the two equations ($\nabla_{a}\nabla_{b}\sigma+P_{ab}\sigma=0$ and $\text{tf}(\nabla_{a}\mu^{b})=0$). In this chapter, we will discuss the settings in which prolongation can be utilised, and thus which equations can be studied by this method. The metrisability equation is prolonged in section 3.4 and we discuss how duality links in with this.

\section{The Projective Tractor Connection}
In the previous chapter we claimed that the two tractor connections given were dual. Here we will prove this.\\\\
\textbf{Claim:}
  \begin{equation*}
    \nabla^{\mathcal{T}}_{a}\begin{pmatrix}
    \mu^{b}\\t
    \end{pmatrix}=\begin{pmatrix}
    \nabla_{a}\mu^{b}+t\delta_{a}^{b}\\
    \nabla_{a}t-P_{ab}\mu^{b}
    \end{pmatrix}
\quad \text{and} \quad
    \nabla_{a}^{\mathcal{T}\ast}\begin{pmatrix}
    \rho\\\nu_{b}
    \end{pmatrix}=\begin{pmatrix}
    \nabla_{a}\rho-\nu_{a}\\
    \nabla_{a}\nu_{b}+P_{ab}\rho
    \end{pmatrix}
\end{equation*}
are dual connections.\\\\
\textbf{Proof:} \\
To show this, we need to show that they satisfy the Leibniz rule
\begin{equation}\label{leib}
    (\nabla_{a}^{\mathcal{T}\ast}U_{B})V^{B}+U_{B}(\nabla_{a}^{\mathcal{T}}V^{B})=\nabla_{a}^{\mathcal{T}}(U_{B}V^{B})
\end{equation}
as this verifies that $\nabla_{a}^{\mathcal{T}\ast}$ is the dual connection to $\nabla_{a}^{\mathcal{T}}$ defined by the Leibniz rule. $U_{B}V^{B}$ is the natural pairing of an element in the tractor bundle with an element in the cotractor bundle. Therefore, $\nabla_{a}^{\mathcal{T}}(U_{B}V^{B})$ is just the exterior derivative of $U_{B}V^{B}$. We introduce notation so that contraction of elements in the tractor bundle and its dual appears as matrix multiplication. Thus $U_{B}V^{B}=\begin{psmallmatrix} \nu_{b}&\rho\end{psmallmatrix}\begin{psmallmatrix}\mu^{b}\\t
\end{psmallmatrix}=\rho t+\nu_{b}\mu^{b} $ where the elements of the dual bundle are written as row vectors in reversed order, for example $\begin{psmallmatrix}
    \rho\\\nu_{b}
    \end{psmallmatrix}$ becomes $\begin{psmallmatrix}
    \nu_{b}&\rho
    \end{psmallmatrix}$.

Now we begin with the left hand side of \ref{leib} and substitute how the connections act on the bundles
\begin{equation*}
\begin{split}
(\nabla_{a}^{\mathcal{T}\ast}U_{B})V^{B}+U_{B}(\nabla_{a}^{\mathcal{T}}V^{B})&=(\nabla_{a}^{\mathcal{T}\ast}\begin{pmatrix}\nu_{b}&\rho \end{pmatrix}) \begin{pmatrix}\mu^{b}\\t
\end{pmatrix}+\begin{pmatrix}\nu_{b}&\rho \end{pmatrix}\nabla_{a}^{\mathcal{T}}\begin{pmatrix}\mu^{b}\\t\end{pmatrix}\\
& =\begin{pmatrix}\nabla_{a}\nu_{b}+P_{ab}\rho&\nabla_{a}\rho-\nu_{a}\end{pmatrix}\begin{pmatrix}\mu^{b}\\t
\end{pmatrix}+\begin{pmatrix}\nu_{b}&\rho \end{pmatrix}\begin{pmatrix}
    \nabla_{a}\mu^{b}+t\delta_{a}^{b}\\
    \nabla_{a}t-P_{ab}\mu^{b}
    \end{pmatrix}\\
    &=\mu^{b}\nabla_{a}\nu_{b}+\mu^{b}P_{ab}\rho+t\nabla_{a}\rho-t\nu_{a}+\nu_{b}\nabla_{a}\mu^{b}+t\nu_{a}+\rho\nabla_{a}t-\rho P_{ab}\mu^{b}\\
    &=\mu^{b}\nabla_{a}\nu_{b}+\nu_{b}\nabla_{a}\mu^{b}+t\nabla_{a}\rho+\rho\nabla_{a}t\\
    &=d_{a}(\rho t+\nu_{b}\mu^{b})
\end{split}
\end{equation*}\\
Thus, the Leibniz rule holds and so \ref{projtrac} and \ref{dualprojtrac} are dual connections.\qed

This is provided as an example that connections can be verified to be dual by the Leibniz rule. We will later see that how a connection acts on the dual bundle can be directly calculated using the Leibniz rule (which we do in section 3.5).

\section{Solutions to Prolonged Equations}

Since parallel sections of the tractor connection are in one-to-one correspondence with solutions to the prolonged equation, there are at most $n+1$ linearly independent solutions to the equation (and hence also for the equation that was prolonged to obtain the cotractor connection) as rank$(\mathcal{E}_{A})=\text{rank}(\mathcal{E}_{a}(-1))+\text{rank}(\mathcal{E}(-1))=n+1$ and of course, the dual bundle has the same rank. In other words, the dimension of the solution space is at most $n+1$. In the flat case, no solutions are obstructed by curvature and therefore the number of linearly independent solutions is equal to the bound. The system can be directly solved by integrating to get one solution and then obtaining the rest through parallel transport of the first one. This means that in a flat setting two equations that give dual connections through prolongation have the same number of solutions.

For any connection obtained through the prolongation of an equation, parallel sections of the connection are in one-to-one correspondence with solutions to the equation. Therefore, the rank of the prolonged system gives an upper bound on the dimension of the solution space, this bound is not necessarily met as curvature can obstruct some of the solutions.

\section{When is Prolongation Successful}
Prolongation is a method originally used for ordinary differential equations. It is a common tool which allows an order $n$ ordinary differential equation to be reduced to a system of $n$ linear first order ordinary differential equations. Partial differential equations however do not directly reduce and it is necessary to introduce new variables for higher derivatives in order to form a closed system. A closed system is the aim of prolongation, we will define this here.

\begin{definition}[Closed System]
A closed system is a system of differential equations for which the first derivatives of each variable are algebraically expressed in terms of the other variables.
\end{definition}

\noindent We have been concerned with projectively-invariant differential equations on manifolds, however as we will discuss in section 3.7 this can be done in other settings. For an equation to be prolonged successfully and give a closed system it must be of finite type. An equation of finite type is one in which the system closes after the introduction of finitely many new variables.

\section{The Case of the Metrisability Equation}
The metrisability equation is
$$tf(\nabla_{a}t^{bc})=0$$ where the rank $2$ contravariant tensor, $t^{bc}$, is symmetric. Non-degenerate solutions to this equation are equivalent to metrics whose Levi-Civita connection lies in the projective class. We can see that this equation gives $n^{3}$ scalar equations, however we can calculate the number of independent components of $t^{bc}$ to be $\frac{n^{2}(n+1)}{2}-n=\frac{n(n-1)(n+2)}{2}$ which is less that $n^{3}$ for $n\geq 3$. Therefore, the system has more scalar equations than individual scalar components so is overdetermined when $n\geq 3$.

Recall that the Weyl tensor contracted with $t^{bc}$ is given by $W_{ab}{}^{c}{}_{d}t^{bd}=R_{ab}{}^{c}{}_{d}t^{bd}-\delta_{a}^{c}P_{bd}t^{bd}+P_{ad}t^{cd}$. Since $\nabla_{a}t^{bc}-\text{trace-part}(\nabla_{a}t^{bc})=0$, we begin by formulating an expression for the trace part of $\nabla_{a}t^{bc}$. Because $t^{bc}$ is symmetric, this will take the form $k\delta_{a}^{b}\nabla_{d}t^{dc}+k\delta_{a}^{c}\nabla_{d}t^{bd}$ for a suitable $k$ which can be verified to be $k=\frac{1}{n+1}$ (this is calculated in \ref{appendixA}). Therefore,
$$\nabla_{a}t^{bc}-\frac{1}{n+1}\delta_{b}^{a}\nabla_{d}t^{dc}-\frac{1}{n+1}\delta_{c}^{a}\nabla_{d}t^{db}=0.$$
Defining $\nu^{c}=-\frac{1}{n+1}\nabla_{d}t^{dc}$ allows us to rewrite the metrisability equation as the following,

\begin{equation}\label{1m}
    \nabla_{a}t^{bc}+\delta_{a}^{b}\nu^{c}+\delta_{a}^{c}\nu^{b}=0.
\end{equation}
Now we take the derivative, $\nabla_{e}$,

$$\nabla_{e}\nabla_{a}t^{bc}+\nabla_{e}\delta_{a}^{b}\nu^{c}+\nabla_{e}\delta_{a}^{c}\nu^{b}=0.$$

\noindent Now with our previous definition of curvature we can rewrite $\nabla_{e}\nabla_{a}t^{bc}$ such that it is in terms of $\nabla_{a}\nabla_{e}t^{bc}$ and a curvature term, hence commuting the covariant derivatives.

$$\nabla_{a}\nabla_{e}t^{bc}+R_{ea}{}^{b}{}_{d}t^{dc}+R_{ea}{}^{c}{}_{d}t^{bd}+\nabla_{e}\delta_{a}^{b}\nu^{c}+\nabla_{e}\delta_{a}^{c}\nu^{b}=0$$

\noindent Then we trace over $e$ and $b$ as we would like to form an equation for $\nabla_{a}\nu^{c}$

$$\nabla_{a}\nabla_{b}t^{bc}+R_{ad}t^{cd}+R_{ba}{}^{c}{}_{d}t^{bd}+\nabla_{a}\nu^{c}+\delta_{a}^{c}\nabla_{b}\nu^{b}=0$$

\noindent and using the definition of $\nu^{c}$ and introducing the Weyl tensor gives
\begin{equation*}-(n+1)\nabla_{a}\nu^{c}+R_{ad}t^{cd}-W_{ab}{}^{c}{}_{d}t^{bd}+P_{ab}t^{bc}-\delta_{a}^{c}P_{bd}t^{bd}+\nabla_{a}\nu^{c}+\delta^{c}_{a}\nabla_{b}\nu^{b}=0\end{equation*}

\noindent We now divide by the coefficient of $\nabla_{a}\nu^{c}$

$$-\nabla_{a}\nu^{c}+\frac{1}{n}R_{ad}t^{cd}-\frac{1}{n}W_{ab}{}^{c}{}_{d}t^{bd}+\frac{1}{n}P_{ab}t^{bc}-\frac{1}{n}\delta_{a}^{c}P_{bd}t^{bd}+\frac{1}{n}\delta^{c}_{a}\nabla_{b}\nu^{b}=0$$

\noindent Then we define $\rho=-\frac{1}{n}\nabla_{a}\nu^{a}+\frac{1}{n}P_{de}t^{ed}$ and substitute it in

$$-\nabla_{a}\nu^{c}+\frac{1}{n}R_{ad}t^{cd}-\frac{1}{n}W_{ab}{}^{c}{}_{d}t^{bd}+\frac{1}{n}P_{ab}t^{bc}-\delta_{a}^{c}\rho=0$$

\noindent And $\frac{1}{n}R_{ab}t^{bc}+\frac{1}{n}P_{ab}t^{bc}=\frac{n-1}{n}P_{ab}t^{bc}+\frac{1}{n}P_{ab}t^{bc}=P_{ab}t^{bc}$. Therefore, we have,
\begin{equation}\label{2m}
    \nabla_{a}\nu^{c}+\delta^{c}_{a}\rho-P_{ab}t^{cb}+\tfrac{1}{n}W_{ab}{}^{c}{}_{d}t^{bd}=0
\end{equation}

\noindent Now we want an equation for $\nabla_{a}\rho$ so we take the covariant derivative $\nabla_{e}$ and using the Leibniz rule on the last three terms ($\nabla_{e}\delta_{a}^{c}=0$).
\begin{equation*}
    \nabla_{e}\nabla_{a}\nu^{c}+\delta^{c}_{a}\nabla_{e}\rho-t^{cb}\nabla_{e}P_{ab}-P_{ab}\nabla_{e}t^{cb}+\tfrac{1}{n}t^{bd}\nabla_{e}W_{ab}{}^{c}{}_{d}+\tfrac{1}{n}W_{ab}{}^{c}{}_{d}\nabla_{e}t^{bd}=0
\end{equation*}

\noindent Again, we swap the covariant derivatives by introducing a curvature term
\begin{equation*}
    \nabla_{a}\nabla_{e}\nu^{c}+R_{ea}{}^{c}{}_{d}\nu^{d}+\delta^{c}_{a}\nabla_{e}\rho-t^{cb}\nabla_{e}P_{ab}-P_{ab}\nabla_{e}t^{cb}+\tfrac{1}{n}t^{bd}\nabla_{e}W_{ab}{}^{c}{}_{d}+\tfrac{1}{n}W_{ab}{}^{c}{}_{d}\nabla_{e}t^{bd}=0
\end{equation*}

\noindent Now let us trace over $c$ and $e$,
\begin{equation*}
    \nabla_{a}\nabla_{c}\nu^{c}+(n-1)P_{ad}\nu^{d}+\nabla_{a}\rho-t^{cb}\nabla_{c}P_{ab}+P_{ab}(n+1)\nu^{b}+\tfrac{n-2}{n}t^{bd}C_{abd}-\tfrac{1}{n}W_{ab}{}^{c}{}_{d}(\delta_{c}^{b}\nu^{d}+\delta_{c}^{d}\nu^{b})=0
\end{equation*}

\noindent Next we want to subtract $\nabla_{a}P_{de}t^{de}$ from both sides in order to substitute in $\rho$.
\begin{equation*}
    \nabla_{a}\nabla_{c}\nu^{c}+(n-1)P_{ad}\nu^{d}+\nabla_{a}\rho-t^{cb}\nabla_{c}P_{ab}+P_{ab}(n+1)\nu^{b}+\tfrac{n-2}{n}t^{bd}C_{abd}-\nabla_{a}(P_{de}t^{de})+\nabla_{a}(P_{de}t^{de})=0
\end{equation*}

\noindent Now we substitute in $\rho$.
\begin{equation*}
    (1-n)\nabla_{a}\rho+(2n-2)P_{ad}\nu^{d}+\tfrac{n-2}{n}t^{bd}C_{abd}+t^{cb}\nabla_{a}P_{cb}-t^{cb}\nabla_{c}P_{ab}=0
\end{equation*}

\noindent And finally after simplifying we have

\begin{equation}\label{3m}
    \nabla_{a}\rho-2P_{ad}\nu^{d}-\tfrac{2}{n}t^{bd}C_{abd}=0.
\end{equation}

\noindent \ref{1m}, \ref{2m}, and \ref{3m} form a closed system and thus the prolongation is complete. The three equations define a connection acting on $S^{2}\mathcal{T}$ with a bijection between solutions of the metrisability equation and parallel sections of the connection.

\begin{equation}\label{prolongationconnection}
D_{a}\begin{pmatrix}
t^{bc}\\\nu^{c}\\\rho
\end{pmatrix}=\begin{pmatrix}
\nabla_{a}t^{bc}+\delta_{a}^{b}\nu^{c}+\delta_{a}^{c}\nu^{b}\\
\nabla_{a}\nu^{c}+\delta^{c}_{a}\rho-P_{ab}t^{cb}+\tfrac{1}{n}W_{ab}{}^{c}{}_{d}t^{bd}\\
\nabla_{a}\rho-2P_{ad}\nu^{d}-\tfrac{2}{n}t^{bd}C_{abd}
\end{pmatrix}
\end{equation}

\noindent \textbf{Remark:} The prolongation of $tf(\nabla_{a}\beta^{bc})=0$ where $\beta^{bc}$ is skew symmetric is done in the flat case in \ref{appendixC}, the reader is encouraged to do this as an exercise before consulting the appendix. 
\section{Calculating the Dual Connection}

Next let us explicitly calculate the dual to \ref{prolongationconnection} using the Leibniz rule. This will give us a connection acting on $S^{2}\mathcal{T}^{\ast}$. Furthermore, if we can find which equation is prolonged to give it then we have also found a link between these two equations. We begin by plugging \ref{prolongationconnection}
into the Leibniz rule.
\begin{equation*}
    \begin{split}
        (\nabla_{a}^{\ast} \begin{pmatrix}
        \sigma & \mu_{c} & \beta_{bc})
        \end{pmatrix}\begin{pmatrix}
        \rho\\ \nu^{c} \\ t^{bc}\end{pmatrix}&=\nabla_{a}(t^{bc}\beta_{bc}+\mu_{c}\nu^{c}+\rho t)-\begin{pmatrix}\sigma & \mu_{c} & \beta_{bc}\end{pmatrix}\nabla_{a}\begin{pmatrix}t^{bc}\\\nu^{c}\\\rho\end{pmatrix}\\
        &=t^{bc}\nabla_{a}\beta_{bc}+\beta_{bc}\nabla_{a}t^{bc}+\mu_{c}\nabla_{a}\nu^{c}+\nu^{c}\nabla_{a}\mu_{c}+\rho\nabla_{a}t+t\nabla_{a}\rho\\ 
        &\quad -\begin{pmatrix}\sigma&\mu_{c}&\beta_{bc}\end{pmatrix}\begin{pmatrix}\nabla_{a}t^{bc}+\delta_{a}^{b}\nu^{c}+\delta_{a}^{c}\nu^{b}\\\nabla_{a}\nu^{c}+\delta^{c}_{a}\rho-P_{ab}t^{cb}+\tfrac{1}{n}W_{ab}{}^{c}{}_{d}t^{bd}\\\nabla_{a}\rho-2P_{ad}\nu^{d}-\tfrac{2}{n}t^{bd}C_{abd}\end{pmatrix}\\
        &=t^{bc}\nabla_{a}\beta_{bc}+\nu^{c}\nabla_{a}\mu_{c}+\rho\nabla_{a}\sigma-2\beta_{ac}\nu^{c} -\mu_{a} \rho\\
        &\quad+\mu_{c}P_{ab}t^{cb}-\mu_{c}\tfrac{1}{n}W_{ab}{}^{c}{}_{d}t^{bd}+\sigma2P_{ad}\nu^{d}+\sigma\tfrac{2}{n}t^{bd}C_{abd}
    \end{split}
\end{equation*}

\noindent Now we want to separate the terms into the terms multiplied by $\rho$, the terms multiplied by $\nu^{c}$, and the terms multiplied by $t^{bc}$. How we have defined multiplication of elements of $S^{2}\mathcal{T}$ and $S^{2}\mathcal{T}^{\ast}$ means that the terms multiplied by $t^{bc}$ will become the top line.  

\begin{equation*}
    D_{a}\begin{pmatrix}
    \beta_{bc}\\\mu_{c}\\\sigma
    \end{pmatrix}=\begin{pmatrix}
    \nabla_{a}\beta_{bc}+\mu_{c}P_{ab}+\mu_{b}P_{ac}-\mu_{e}W_{ab}{}^{e}{}_{c}t^{bc}-\mu_{e}W_{ac}{}^{e}{}_{b}t^{bc}+\sigma\frac{2}{n}C_{abc}+\sigma\frac{2}{n}C_{acb}\\
    \nabla_{a}\mu_{c}-2\beta_{ac}+2P_{ac}\sigma\\
    \nabla_{a}\sigma-\mu_{a}
    \end{pmatrix}.
\end{equation*}

\noindent Now we have a connection defined on $S^{2}\mathcal{T}^{\ast}$.

\section{Using the Tractor Connection to get Further Prolongations}
The above example illustrates that prolongation can quickly become a tedious and laborious task. It turns out that previously prolonged systems can be used to simplify the prolongation of higher order equations. This is what we shall study in this section.

The tractor bundle whose elements we write conveniently in our notation of $\begin{pmatrix}
\rho\\\nu^{b}
\end{pmatrix}$ is formally represented by canonical bundle maps $X^{A}$, $Y_{A}$, $Z_{A}^{a}$, and $W^{A}_{a}$.

\begin{equation*}
    \begin{split}
    \nabla_{e}\begin{pmatrix}
t^{bc}\\\nu^{c}\\\rho
\end{pmatrix}&=\nabla_{e}(W_{a}^{B}W_{d}^{C}t^{ad})+\nabla_{e}(X^{B}W^{C}_{a}\nu^{a})+\nabla_{e}(X^{C}W^{B}_{a}\nu^{a})+\nabla_{e}(X^{B}X^{C}\rho)\\
        &=\nabla_{e}(W_{a}^{B}W_{d}^{C}t^{ad}+X^{B}W^{C}_{a}\nu^{a}+X^{C}W^{B}_{a}\nu^{a}+X^{B}X^{C}\rho)\\&=t^{ad}W_{a}^{B}\nabla_{e}W_{d}^{C}+W_{d}^{C}t^{ad}\nabla_{e}W_{a}^{B}+W_{a}^{B}W_{d}^{C}\nabla_{e}t^{ad}\\&\quad +X^{B}\nu^{a}\nabla_{e}W^{C}_{a}+W^{C}_{a}\nu^{a}\nabla_{e}X^{B}+X^{B}W^{C}_{a}\nabla_{e}\nu^{a}\\&\quad+\nu^{a}X^{C}\nabla_{e}W^{B}_{a}+X^{C}W^{B}_{a}\nabla_{e}\nu^{a}+W^{B}_{a}\nu^{a}\nabla_{e}X^{C}\\&\quad+\rho X^{B}\nabla_{e}X^{C}+X^{C}\rho\nabla_{e}X^{B}+X^{B}X^{C}\nabla_{e}\rho
    \end{split}
\end{equation*}

Now we can substitute in $\nabla_{e}W^{B}_{b}$ and $\nabla_{e}X^{B}$. This is how the tractor connection on the tractor bundle is defined and we will calculate it from how the tractor connection acts on the tractor bundle.

Elements of the tractor bundle are written formally as $W^{B}_{b}\mu^{b}+X^{B}t$. And we have that $\nabla_{a}^{\mathcal{T}}(W^{B}_{b}\mu^{b}+X^{B}t)=(\nabla_{a}\mu^{b}+t\delta_{a}^{b})W^{B}_{b}+(\nabla_{a}t-P_{ab}\mu^{b})X^{B}=W^{B}_{b}\nabla_{a}\mu^{b}+tW_{a}^{B}+X^{B}\nabla_{a}t-P_{ab}\mu^{b}X^{B}$. Thus we identify the relevant bits and get that $\nabla_{a}W^{B}_{b}=-P_{ab}X^{B}$ and $\nabla_{a}X^{B}=W_{a}^{B}$

Now we will substitute this in.

\begin{equation*}
    \begin{split}
    \nabla_{e}\begin{pmatrix}
t^{bc}\\\nu^{c}\\\rho
\end{pmatrix}&=\nabla_{e}(W_{a}^{B}W_{d}^{C}t^{ad})+\nabla_{e}(X^{B}W^{C}_{a}\nu^{a})+\nabla_{e}(X^{C}W^{B}_{a}\nu^{a})+\nabla_{e}(X^{B}X^{C}\rho)\\
        &=t^{ad}W_{a}^{B}(-P_{ed}X^{C})+W_{d}^{C}t^{ad}(-P_{ea}X^{B})+W_{a}^{B}W_{d}^{C}\nabla_{e}t^{ad}\\&\quad +X^{B}\nu^{a}(-P_{ea}X^{C})+W^{C}_{a}\nu^{a}W_{e}^{B}+X^{B}W^{C}_{a}\nabla_{e}\nu^{a}\\&\quad+\nu^{a}X^{C}(-P_{ea}X^{B})+X^{C}W^{B}_{a}\nabla_{e}\nu^{a}+W^{B}_{a}\nu^{a}W_{e}^{C}\\&\quad+\rho X^{B}W_{e}^{C}+X^{C}\rho W_{e}^{B}+X^{B}X^{C}\nabla_{e}\rho\\
        &=W^{B}_{a}W^{C}_{d}(\nabla_{e}t^{ad}+\delta_{e}^{a}\nu^{d}+\delta_{e}^{d}\nu^{a})+W^{C}_{d}X^{B}(\nabla_{e}\nu^{d}-P_{ea}t^{ad}+\delta_{e}^{d}\rho)\\&\quad +W^{B}_{a}X^{C}(-P_{ed}t^{ad}+\nabla_{e}\nu^{d}+\delta_{e}^{d}\rho)+X^{B}X^{C}(\nabla_{e}\rho-P_{ea}\nu^{a}-P_{ea}\nu^{a})
    \end{split}
\end{equation*}

\noindent Therefore,

\begin{equation}\label{3.6}
    \nabla_{e}\begin{pmatrix}
t^{bc}\\\nu^{c}\\\rho
\end{pmatrix}=\begin{pmatrix}
\nabla_{e}t^{bc}+\delta_{e}^{b}\nu^{c}+\delta_{e}^{c}\nu^{b}\\
\nabla_{e}\nu^{c}-P_{eb}t^{bc}+\delta_{e}^{c}\rho\\
\nabla_{e}\rho-2P_{eb}\nu^{b}

\end{pmatrix}
\end{equation}

\noindent This is how the tractor connection acts on $S^{2}\mathcal{T}$. This is a distinct connection to the one we prolonged in section 3.3, however both act on $S^{2}\mathcal{T}$. It can be illuminating then to study the conditions under which they agree, as this will tell us about the geometry. Solutions to the metrisability equation are thus given by,

\begin{equation}\label{3.7}
    \nabla_{a}^{\mathcal{T}}\begin{pmatrix}
    t^{bc}\\
    \nu^{c}\\
    \rho\\
    \end{pmatrix}=\begin{pmatrix}0\\
    -\frac{1}{n}W_{ab}{}^{c}{}_{d}t^{bd}\\
    \frac{2}{n}C_{abd}t^{bd}
    \end{pmatrix}
\end{equation}

\noindent The two connections \ref{prolongationconnection} and \ref{3.6} then agree when $\begin{pmatrix}
0\\-\frac{1}{n}W_{ab}{}^{c}{}_{d}t^{bd}\\
\frac{2}{n}C_{abd}t^{bd}
\end{pmatrix}$ vanishes. Clearly this occurs in a flat setting. When we say $t^{bc}$ is equivalent to a metric, we mean that $g^{bc}=\text{sgn}(\tau)\tau t^{bc}$ is the inverse of a metric for $\tau$, a suitable determinant of $t^{bc}$ \cite{FloodGoverMetrisability}. It is easy to show that this occurs if and only if $t^{bc}$ is equivalent to an Einstein metric. We provide a proof. Since $t^{bc}$ is a solution to the metrisability equation, it is equivalent to a metric, $g^{bc}$, and therefore we can deal with $g^{bc}$ instead of $t^{bc}$ in our proof, and we will work with the Levi-Civita connection of $g$, $\nabla^{g}$, since it is in the projective class.

\noindent \textbf{Proof:} Suppose $t^{bc}$ is a solution to the metrisability equation.

\begin{equation*}
\begin{split}
    W_{ab}{}^{c}{}_{d}g^{bd}&=R_{ab}{}^{c}{}_{d}g^{bd}-\delta_{a}^{c}P_{bd}g^{bd}+P_{ad}g^{cd}\\
    &=R_{abed}g^{bd}g^{ce}-\frac{1}{n-1}\delta_{a}^{c}R_{bd}g^{bd}+\frac{1}{n-1}R_{ad}g^{cd}\\
    &=R_{bade}g^{bd}g^{ce}-\frac{1}{n-1}Rg_{ad}g^{cd}+\frac{1}{n-1}R_{ad}g^{cd}\\
    &=R_{ae}g^{ce}-\frac{1}{n-1}Rg_{ad}g^{cd}+\frac{1}{n-1}R_{ae}g^{ce}\\
    &=\frac{n}{n-1}g^{cd}(R_{ad}-\frac{1}{n}g_{ad}R)
\end{split}
\end{equation*}

\noindent Where we have used the symmetries of the Riemmanian curvature, and used the metric to raise and lower indices. Now suppose $g$ is Einstein. Then $R_{ab}=\lambda g_{ab}$. Taking the trace yields $g^{ab}R_{ab}=\lambda g_{ab}g^{ab}$ which means that $R=n\lambda$, and so $R_{ab}=\frac{1}{n}Rg_{ab}$. Thus, $W_{ab}{}^{c}{}_{d}g^{bd}=\frac{n}{n-1}g^{cd}(R_{ad}-\frac{1}{n}g_{ad}R)=0$.

Now to prove that the metric is Einstein if the right hand side of \ref{3.7} vanishes, suppose $W_{ab}{}^{c}{}_{d}g^{bd}=0$ and $C_{abd}g^{bd}=g^{bd}\nabla_{c}W_{ab}{}^{c}{}_{d}=0$. This system is fully determined by $W_{ab}{}^{c}{}_{d}g^{bd}=0$. Therefore, $\frac{n}{n-1}g^{cd}(R_{ad}-\frac{1}{n}g_{ad}R)=0$. And so $R_{ad}=\frac{1}{n}g_{ad}R$ which means that $R_{ad}=\lambda g_{ad}$. Thus, $g_{ab}$ is Einstein. Therefore, the right hand side of \ref{3.7} vanishes if and only if $g_{ab}$ is Einstein.

It only remains to show that the Cotton tensor vanishes if $g$ is Einstein. Supposing $g$ is Einstein, we have that $R_{ad}=\lambda g_{ad}$. Thus $P_{ad}=\frac{1}{n-1}\lambda g_{ad}$ which means that $\nabla_{b}P_{ad}=\frac{1}{n-1}\lambda\nabla_{b}g_{ad}=0$ since the Levi-Civita connection preserves the metric. \qed

\section{Prolongation in other settings}
Prolongation as a method to study the underlying geometry of a space and an equation is not limited to the context of  projective geometry. Prolongation has been studied for conformal geometry and is also used to obtain the conformal tractor connection from a conformally invariant equation. This can be seen in \cite{CurryGoverGR} and the reader is directed to \cite{EastwoodCDGNotes} for an introduction to conformal differential geometry.

\chapter*{Conclusion} In this dissertation we discuss the links between projective tractor calculus and prolongation. We have shown that the projective tractor connection and prolongation of the metrisability equation provide two distinct connections on $S^{2}\mathcal{T}$. We prove they agree if and only if there is an Einstein Levi-Civita connection in the projective class of the projective manifold. This is an important result as projective classes which include an Einstein metric have a special role in projective differential geometry. This dissertation serves as a resource that gives examples of prolongations step by step and also introduces differential geometry and tractor calculus. The role of prolongation in these results shows that it can be an illuminating method to study equations and infer things about the underlying geometries. This dissertation bridges some gaps between introductory differential geometry and topics that are under development in research papers currently (prolongation and projective geometry).

\addcontentsline{toc}{chapter}{Conclusion}
\markboth{CONCLUSION}{CONCLUSION}

\nocite{hawking_ellis_1973}
\nocite{shlomo}
\printbibliography

\appendix
\chapter{Calculating $\boldsymbol{k}$}\label{appendixA}
We begin by formulating an expression for the trace part of $\nabla_{a}t^{bc}$ which will take the form $k\delta_{a}^{b}\nabla_{d}t^{dc}+k\delta_{a}^{c}\nabla_{d}t^{bd}$ for a suitable $k$. Hence,

$$\nabla_{a}t^{bc}-k\delta_{a}^{b}\nabla_{d}t^{dc}-k\delta_{a}^{c}\nabla_{d}t^{bd}=0$$

\noindent Now,  we take the trace over $a$ and $b$.

$$\delta_{b}^{a}\nabla_{a}t^{bc}-\delta_{b}^{a}k\delta_{a}^{b}\nabla_{d}t^{dc}-\delta_{b}^{a}k\delta_{a}^{c}\nabla_{d}t^{bd}=0$$

\noindent And then since $\delta_{a}^{b}\delta_{b}^{a}=n$,

$$\delta_{b}^{a}\nabla_{a}t^{bc}-nk\nabla_{d}t^{dc}-k\delta_{b}^{c}\nabla_{d}t^{bd}=0$$
Now we simplify,

$$\nabla_{b}t^{bc}-nk\nabla_{d}t^{dc}-k\nabla_{d}t^{cd}=0$$
and utilise the symmetry of $t^{cd}$,
$$\nabla_{d}t^{dc}-nk\nabla_{d}t^{dc}-k\nabla_{d}t^{dc}=0$$.
Finally, this allows us to form an expression for $k$,

$$(1-nk-k)\nabla_{d}t^{dc}=0$$

$$\implies 1-nk-k=0$$

$$\implies k=\frac{1}{n+1}$$

\chapter{Proof of the Koszul formula}\label{appendixB}
Suppose we have a connection with the following properties

\begin{enumerate}
    \item Torsion freeness: $[X,Y]=\nabla_{X}Y-\nabla_{Y}X$
    \item Compatibility with the metric: $Xg(Y,Z)=g(\nabla_{X}Y,Z)+g(Y,\nabla_{X}Z).$ 
\end{enumerate}
For such a connection we compute $Xg(Y,Z)+Yg(X,Z)-Zg(X,Y)$ by substituting in the metric compatibility condition. This gives
$$g(\nabla_{X}Y,Z)+g(Y,\nabla_{X}Z)+g(\nabla_{Y}Z,X)+g(\nabla_{Y}X,Z)-g(\nabla_{Z}X,Y)-g(X,\nabla_{Z}Y)$$
Next we rearrange using linearity of the metric.
$$g(\nabla_{X}Y+\nabla_{Y}X,Z)+g(\nabla_{X}Z-\nabla_{Z}X,Y)+g(X,\nabla_{Y}Z-\nabla_{Z}Y)$$
Substituting in the torsion free condition gives

$$g(\nabla_{X}Y+\nabla_{Y}X+(\nabla_{X}Y-\nabla_{X}Y),Z)+g([X,Z],Y)+g([Y,Z],X)$$
Finally, rearranging the first term gives
$$2g(\nabla_{X}Y,Z)+g([X,Y],Z)+g([X,Z],Y)+g([Y,Z],X)$$
By equating this to the expression we computed we get 
\begin{equation*}
\begin{split}
    g(\nabla_{X}Y,Z)=\tfrac{1}{2}\big\{&Xg(Y,Z)+Yg(X,Z)-Zg(X,Y)-g([X,Y],Z)\\&-g([Z,X],Y)-g([Y,Z],X)\big\}
\end{split}
\end{equation*}
which is the Koszul formula. This defines a connection $\nabla:(X,Y)\to \nabla_{X}Y$ whic is torsion free and compatible with the metric. Now, suppose $g(\nabla_{X}Y,Z)=g(W,Z)$ for arbitrary $Z$, then $W=\nabla_{X}Y$. This shows that the two conditions - torsion freeness and metric compatibility - guarantee existence and uniqueness of this connection.
\chapter{Extra Prolongation Example}\label{appendixC}
Recall that a flat connection is one whose associated curvature vanishes everywhere on the manifold. On a manifold, $M$, with flat connection, $\nabla$, consider the equation 
$$tf(\nabla_{a}\beta^{bc})=0$$
where the rank $2$ contravariant tensor, $\beta^{bc}$, is skew-symmetric. Since $\nabla_{a}\beta^{bc}-\text{tr}(\nabla_{a}\beta^{bc})=0$, we begin by formulating an expression for the trace of $\nabla_{a}\beta^{bc}$ which will take the form $k\delta_{a}^{b}\nabla_{d}\beta^{dc}+k\delta_{a}^{c}\nabla_{d}\beta^{bd}$ for a suitable $k$. Hence,

$$\nabla_{a}\beta^{bc}-k\delta_{a}^{b}\nabla_{d}\beta^{dc}-k\delta_{a}^{c}\nabla_{d}\beta^{bd}=0$$

Now,  we take the trace over $a$ and $b$.

$$\delta_{b}^{a}\nabla_{a}\beta^{bc}-\delta_{b}^{a}k\delta_{a}^{b}\nabla_{d}\beta^{dc}-\delta_{b}^{a}k\delta_{a}^{c}\nabla_{d}\beta^{bd}=0$$

And then since $\delta_{a}^{b}\delta_{b}^{a}=n$,

$$\delta_{b}^{a}\nabla_{a}\beta^{bc}-nk\nabla_{d}\beta^{dc}-k\delta_{b}^{c}\nabla_{d}\beta^{bd}=0$$
Now we simplify,

$$\nabla_{b}\beta^{bc}-nk\nabla_{d}\beta^{dc}-k\nabla_{d}\beta^{cd}=0$$
and utilise the skew-symmetry of $\beta^{cd}$,
$$\nabla_{d}\beta^{dc}-nk\nabla_{d}\beta^{dc}+k\nabla_{d}\beta^{dc}=0$$.
Finally, this allows us to form an expression for $k$,
$$\implies 1-nk+k=0$$

$$\implies k=\frac{1}{n-1}$$

Therefore,
$$\nabla_{a}\beta^{bc}-\frac{1}{n-1}\delta_{b}^{a}\nabla_{d}\beta^{dc}-\frac{1}{n-1}\delta_{c}^{a}\nabla_{d}\beta^{bd}=0.$$

Hence,
$$\nabla_{a}\beta^{bc}-\frac{1}{n-1}\delta_{b}^{a}\nabla_{d}\beta^{dc}+\frac{1}{n-1}\delta_{c}^{a}\nabla_{d}\beta^{db}=0.$$
Defining $\nu^{c}=\frac{1}{n-1}\nabla_{d}\beta^{dc}$ allows us to rewrite the problem as the following,
$$\nabla_{a}\beta^{bc}-\delta_{a}^{b}\nu^{c}+\delta_{a}^{c}\nu^{b}=0.$$
Taking the derivative and then the trace yields
$$\nabla_{a}\nu^{c}+\frac{1}{n-2}\delta_{a}^{c}\nabla_{b}\nu^{b}=0.$$
Now we can repeat the process of introducing a variable by defining $\rho = \frac{1}{n-2}\nabla_{b}\nu^{b}$. Hence,

$$(n-1)\nabla_{a}\rho=0,$$ which allows us to conclude that $\nabla \rho=0$. Now putting the three equations for the covariant derivatives of $t, \nu$ and $\rho$ 

\begin{equation*}
    D_{a}\begin{pmatrix}t^{bc}\\\nu^{b}\\\rho\end{pmatrix}\coloneqq\begin{pmatrix}\nabla_{a}t^{bc}-\delta_{a}^{b}\nu^{c}+\delta_{a}^{c}\nu^{b}\\\nabla_{a}\nu^{b}-\delta_{a}^{b}\rho\\\nabla_{a}\rho\end{pmatrix}\end{equation*}

\noindent This bundle is a section of $\mathcal{E}^{[ab]} \oplus\mathcal{E}^{c}\oplus\mathcal{E}$.

\end{document}